\documentclass[12pt,oneside,reqno]{article}
\usepackage{amsmath}
  \usepackage{paralist}
  \usepackage{graphics}
  \usepackage{epsfig}
  \usepackage{float}
  \usepackage{graphicx}
 \usepackage{epstopdf}
 \usepackage{amssymb}
 \usepackage{cite}
 \usepackage{mathrsfs}
 \usepackage{amsthm}
 \usepackage{color}
 \usepackage{caption}
 \usepackage{subfigure}
 \usepackage[colorlinks=true]{hyperref}
\hypersetup{urlcolor=blue, citecolor=blue}
\usepackage[utf8]{inputenc}
 \usepackage[margin=1.0in]{geometry}
\numberwithin{equation}{section}

\newtheorem{theorem}{Theorem}[section]
\newtheorem{corollary}{Corollary}
\newtheorem{lemma}[theorem]{Lemma}

\newtheorem{physical conclusion}{Physical Conclusion}

\newtheorem{definition}[theorem]{Definition}
\newtheorem{remark}{Remark}

\title{The well-posedness and the regularity of global attractor for a couple stress fluid through porous layer with the local thermal non-equilibrium effect}
\author{
Liang Li
\thanks{2019322010018@stu.scu.edu.cn }
\quad\quad
Lan Jia
\thanks{Corresponding author: jia8l8l@163.com}
\\ \footnotesize * College of Mathematics,
\\ \footnotesize Sichuan University
\\ \footnotesize Chengdu, Sichuan 610065, China
\\ \footnotesize \dag College of Applied Mathematics,
\\ \footnotesize Chengdu University of Information Technology
\\ \footnotesize Chengdu, Sichuan 610225, China
\thanks{The work was supported by the National Nature Science Foundation of China (11901408) and (11711306).}
}
\medskip

\begin{document}
\date{}
\maketitle
\begin{abstract}
  In the article, we aim to investigate the well-posedness of solution and the regularity of the global attractor for the couple stress fluid in saturated porous media with the local thermal non-equilibrium effect. To be more specific,  we firstly show the existence and uniqueness of global weak solution to the model by making use of the standard Galerkin method. Second, relying on verifying the uniformly compact condition required, we prove the existence of the global attractor of the model in the space where the weak solution resides. Finally, we improve the regularity of global attractor by uniformly compact condition and obtain the $C^{\infty}$ attractor for the model.
\end{abstract}

\section{Introduction}

%Porous media, as is well-known, are ubiquitous in our lives, such as biological tissues (bone, skin), building materials (sand, cement, brick), supplies (ceramics). The motivation for studying porous media is to better understand the properties of such materials and to use them to improve the quality of our lives. Among the subjects related to porous media, global warming is one of current topics and porous media are involved there in connection with ice melting\cite{Carr2003} or carbon dioxide storage\cite{Hill2014}. There are a lot of applications related to porous media, like bone modeling\cite{Eringen2004}, or aircraft or motor car production\cite{DELEGLISE20072317} and so on.
%The study of thermal convection in porous media has important practical significance, such as exploration of geothermal resources, techniques for preventing the spread of pollution sources in underground aquifers, early warning and protection technology of coal seam fire.  Thermal convection in porous media can be simply shown in \autoref{fig1}: the region where the fluid flows is divided into two layers: the one with and the other one without medium, and the bottom of the area is heated. 
The dynamics of fluid through porous media has attracted some of applied mathematicians like Hill and Straughan \cite{Hill2014}\cite{Straughan2005}, and still is a vigorously active research area.
The widely used momentum equation for Newtonian fluid through porous media
is either the Brinkman equation or the Darcy equation derived from the Navier-Stokes equation by  statistical  averages  and  simplifications  of  the  complicated  microscopic flow picture.  Couple stress fluid, developed by V. K. Stokes\cite{Stokes1984} and discussed in detail by himself in his monograph\cite{Stokes1984a}, is a typical non-Newtonian fluid and allows for the polar effects. %such as the presence of couple stresses and body couples. 
The momentum equation \cite{MALASHETTY2009781}\cite{Shivakumara2009} for couple stress fluid in a saturated porous layer reads
%which is more precisely than Newtonian fluid for describing the characteristics of fluid and it has received much attention due to its wide applications in industries.
%The momentum equations of Newtonian fluid through 
% for example, exotic lubricants, the extrusion of polymer fluids, solidification of liquid crystals and so on.
\begin{align}\label{e-c}
\frac{\rho_{0}}{\varepsilon}\frac{\partial \mathbf{v}}{\partial t}+\frac{\mu_{f}}{K}\mathbf{v}
=\frac{\mu_{c}}{K}\Delta \mathbf{v}-\nabla \mathbf{p} -\rho_{0}(1-\beta(T_{f}-T_{u}))\mathbf{g},
\end{align}
where  $\mathbf{v}=(v_{1},v_{2})$ is the velocity field, $T_{f}$ and $T_{u}$ denote the fluid and upper surface temperatures respectively, $\mathbf{p}$ means the pressure field and $\mathbf{g}=(0,-g)$ is the gravitational acceleration. The other numbers in the model are as follows: the porosity $\varepsilon$, the permeability $K$, the viscosity $\mu_{f}$, the coupling stress viscosity $\mu_{c}$, the thermal expand coefficient $\beta$ and  the density $\rho_{0}$ at the initial moment.

The study of thermal convection in porous media has important practical significance, such as exploration of geothermal resources, techniques for preventing the spread of pollution sources in underground aquifers, early warning and protection technology of coal seam fire. However, a lot of work on convection mainly focus on Newtonian fluid, see. e.g. \cite{2006Convection}\cite{Vafai2005}. Now, there is increasingly number of literature considering the convection in couple stress fluid, see. e.g.\cite{2007An}  \cite{Sunil2002}. It is worth noting that the energy models among the above literatures about couple stress fluid are considered under the uniform temperature gradient. However, the local thermal non-equilibrium (LTNE) effect has to be taken into account for establishing the energy equations which is better suitable for the practical situations such as high-speed flows or large temperature differences between the fluid and solid phases. Continuum theories for the LTNE effect on the flow in porous materials appear to have started in the late 1990's, cf. the work of \cite{MINKOWYCZ19993373} where a modified energy equation that can be solved for very early departures from local thermal equilibrium(LTE) conditions is presented while assuming the velocity field is known from the solution of continuity and momentum equation, and \cite{Nield2002}. Thus, for describing the energy equation,  the two-field model \cite{Straughan2015}(page 34)\cite{D.A.Nield2017} 
\begin{align}\label{e-c1}
\begin{cases}
\varepsilon(\rho c)_{f}\frac{\partial T_{f}}{\partial t}+(\rho c)_{f}(\mathbf{v}\cdot \nabla)T_{f}
=\varepsilon \kappa_{f}\Delta T_{f}+h(T_{s}-T_{f}),
\\
(1-\varepsilon)(\rho c)_{s}\frac{\partial T_{s}}{\partial t}
=(1-\varepsilon)\kappa_{s}\Delta T_{s}-h(T_{s}-T_{f})
\end{cases}
\end{align}
is employed, where the energy equations are coupled by the terms that account for the heat lost to or gain form the other phase, $T_{s}$ denotes the solid temperature, $c$ is the specific heat, the heat conductivity $\kappa$ with subscript $f$ and $s$ meaning fluid and solid phase respectively and the inter-phase heat transfer coefficient $h$.

%In current paper, by adding up the acceleration term to the Brinkman's equation\cite{Brinkman1949} and assuming that the coupled effect of temperature on velocity field is linear with respect to temperature,i.e. Boussinesq approximation\cite{Boussinesq1897}, the momentum equation studied is obtained. In contrast with assuming the temperature gradient at any location between the fluid phase and solid phase can be negligible, here 

%Thermal convection in saturated porous media has drawn considerable attention  from a wide range of scientists and engineers, in addition politicians and economists who recognize its importance in practical applications, such as nuclear waste repository, mantle convection and so on. In contrast with assuming the temperature gradient at any location between the fluid phase and solid phase can be negligible, here the local thermal non-equilibrium (LTNE) effect is taken account into establishing the energy equations which is better suitable for the practical situations such as high-speed flows or large temperature differences between the fluid and solid phases. Thus, for describing the energy equation, the two-field model\cite{Straughan2015}(page 34)\cite{D.A.Nield2017}, where the energy equations are coupled by the terms  that account for the heat lost to or gain form the other phase, is employed.

There are many studies involved with couple stress fluid through porous layer using a LTNE model: M.S.Malashetty\cite{Malashetty2009} obtained the condition for the onset of convection by the linear stability theory and also presented asymptotic analysis for different values of the inter-phase heat transfer; Sunli\cite{Sunil2019} showed the equivalence of nonlinear stability threshold and linear instability boundary. Recently, Quan. W\cite{Quangwangjialan2021}studied the stability and transition of the LTNE model of coupled stress fluid. For more information about related model, please refer to the references\cite{Sunil2014}\cite{Kumar2021}. However, there is  lack of work on well-posedness, i.e. the existence and uniqueness of solution, of related model. We will consider the problem later.
 To the best of our knowledge, there is little literature about the existence and regularity of the 
 global attractor for the couple stress fluid in saturated porous layer with the LTNE effect. The study of global attractor\cite{Robinson2009} which is defined as the maximal compact invariant set or the minimal set which uniformly attracts all bounded set can provide deep insight in the long-time behavior of dynamics of the model investigated. In general, there are two ways to prove the existence of the global attractor: Condition(C)\cite{2007} (page 106) \cite{Ma2002} and uniformly compact condition.  Although the C-condition compared with the uniformly compact condition is easier to be verified because there is no need to perform any operations in higher regularity space, it is more convenient to employ the uniformly compact condition to obtain our desire under the help of semigroup. There are a lot of work on the existence of attractor by condition(C), see.e.g. \cite{Zhang_2009}\cite{Ma2002}\cite{Luo2012}.

 Inspired by the work \cite{2011} (page 512) where Ma proved the existence of the $C^{\infty}$-attractor of the 2-dimension Boussinesq equation, this work aim to study the existence of the $C^{\infty}$-attractor for the couple stress fluid in saturated porous layer with the LTNE effect. The approach employed 
in present article is as follows. Firstly, we show the existence of the global weak solution  by using the standard Galerkin method\cite{Evans2010}; 
secondly, relying on verifying the uniformly compact condition required, we prove the existence of the global attractor of the model in the space where the weak solution resides; finally, instead of improving the regularity of the global weak solution by the interpolation theorem\cite{EliasM.Stein2011} used in \cite{2011}, we improve the regularity by estimating the expression of solution directly,  improve the regularity of the global attractor by iterative use of uniformly compact condition
 and obtain the  $C^{\infty}$ attractor for the model. For more literatures on the regularity of the attractor, see. e. g. \cite{Zhang2008}\cite{SONG200953}.

%secondly, the global weak solution is obtained by improving the regularity of the weak solution with respect to the time variable; then the solution can be expressed in form of semigroup. With the help of the expression of the global weak solution, we prove the existence and regularity of global attractor of the model by relying on verifying the uniformly compact condition required.

 %Because the model studied is in 2-dimension space, the uniqueness of weak solution can be proved. Thus, the global solution is also unique.  Global attractor in the space coinciding with the space of weak solution is proved by the uniformly compact condition.  

%However, in contrast with Ma's work where he obtained a sufficient condition for increasingly improving the regularity of global weak solution by the interpolation theorem, here we improved the regularity by estimating the expression of solution directly. Finally, by iterative use of uniformly compact condition, the existence of $C^{\infty}-$attractor, of which the uniqueness comes from the uniqueness of the global attractor, for the model is proved. 

The rest of this article is  arranged as follows. In Section 2, we introduce the mathematical model studied in this paper and its dimensionless form and give some mathematical settings. In section 3, we study the well-posedness of the model solution, that is, we prove the existence of the unique global weak solution and the global solution respectively.
Next, we investigate the existence of the model attractor in detail in Sections 4 and Sections 5. We first prove that there is a global attractor in the weak solution space by using the uniformly compact method. By using the expression of the weak solution and the iterative method, we then improve the regularity of the weak solution and obtain the classical solution of the model under certain conditions. At the same time, we prove that there is a global attractor in the classical solution space.
Finally, we summarize the results of the work in the article in section 6.
% in section 2, we introduce, such as the meaning of mathematical description and so on,  and deal with, like non-dimensionalization and using stream function, the model we consider; in section 3, the well-posedness of solution, including the existence and uniqueness of weak-solution and global solution, has been obtained; in section 4, the existence of attractor in $\mathbf{Y}$ can be proved; in section 5, after iterative improving the regularity of global solution, the $C^{\infty}-$attractor is also obtained by iteration; in section 7, we end the article with a summary of results.

\section{Mathematic setting}
\subsection{Mathematical model}

\begin{figure}[H]
  \centering
  \includegraphics[width=3in]{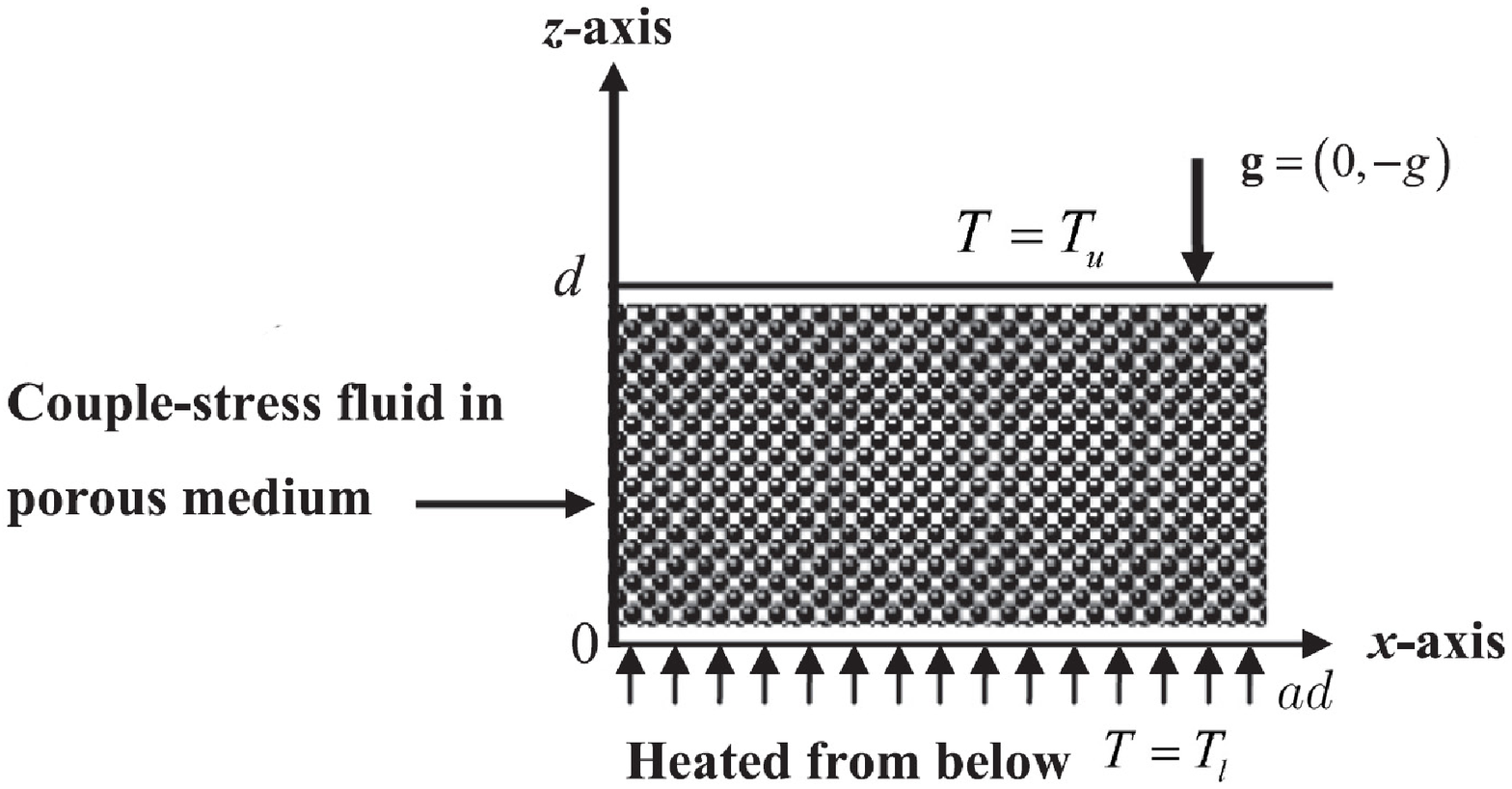}
  \caption{}\label{fig1}
\end{figure}

We consider the 2-D incompressible couple stress fluid model of saturated porous media where the region is a rectangle with depth d and width ad. The fluid is heated from below and cooled from above (see \autoref{fig1}). The temperature of lower surface and of upper surface are held at $T=T_{l}$ and $T=T_{u} (<T_{l})$, respectively.
 %By adding up the acceleration term to the Brinkman's equation which is not just extension of Darcy's Law, the momentum equation we consider can be obtained under some simplific assumptions.
Combining \eqref{e-c}, \eqref{e-c1} and continuity equation, 
the basic governing equations are as follows:
\begin{align}\label{b1}
\begin{cases}
\frac{\rho_{0}}{\varepsilon}\frac{\partial \mathbf{v}}{\partial t}+\frac{\mu_{f}}{K}\mathbf{v}
=\frac{\mu_{c}}{K}\Delta \mathbf{v}-\nabla \mathbf{p} -\rho_{0}(1-\beta(T_{f}-T_{u}))\mathbf{g},
\\
\varepsilon(\rho c)_{f}\frac{\partial T_{f}}{\partial t}+(\rho c)_{f}(\mathbf{v}\cdot \nabla)T_{f}
=\varepsilon \kappa_{f}\Delta T_{f}+h(T_{s}-T_{f}),
\\
(1-\varepsilon)(\rho c)_{s}\frac{\partial T_{s}}{\partial t}
=(1-\varepsilon)\kappa_{s}\Delta T_{s}-h(T_{s}-T_{f}),
\\
 \nabla\cdot \mathbf{v}=0.
\end{cases}
\end{align}

%In the mathematical description, $\mathbf{v}=(v_{1},v_{2})$ is the velocity field, $T_{f}$ and $T_{s}$ denote the fluid and solid temperatures respectively, $\mathbf{p}$ means the pressure field and $\mathbf{g}=(0,-g)$ is the gravitational acceleration. The other numbers in the model are as follows: the porosity $\varepsilon$, the permeability $K$, the viscosity $\mu_{f}$, the coupling stress viscosity $\mu_{c}$, the specific heat $c$, the heat conductivity $\kappa$ with subscript $f$ and $s$ meaning fluid and solid phase respectively, the inter-phase heat transfer coefficient $h$, the thermal expand coefficient $\beta$ and the density $\rho_{0}$ at the initial moment.

Now by the following transformations to make the equation \eqref{b1} dimensionless:
\begin{align}\label{b2}
\begin{aligned}
&(x,z)=(x^\ast,z^\ast),\ \
(v_{1},v_{2})=\frac{\varepsilon \kappa_{f}}{(\rho c)_{f}}(v_{1}^\ast,v_{2}^\ast),\ \
 p =\frac{\kappa_{f}\mu_{f}}{(\rho c)_{f}K} p^\ast ,
 \\
&T_{f}=(T_{l}-T_{u})\theta,\ \
T_{s}=(T_{l}-T_{u})\phi,\ \
t=\frac{(\rho c)_{f}}{\kappa_{f}}t^\ast.
\end{aligned}
\end{align}

After taking the curl of the first dimensionless equation and ignoring the superscript, equation \eqref{b1} is as follows:
\begin{align}\label{b4}
\begin{aligned}
 \frac{Da}{Pr}\frac{\partial \Delta\psi}{\partial t}
&=(C\Delta-1)\Delta\psi+Ra\frac{\partial \theta}{\partial x},
\\
 \frac{\partial \theta}{\partial t}
&=\Delta\theta+\lambda(\phi-\theta)-J(\psi,\theta),
\\
 \alpha\frac{\partial \phi}{\partial t}
&=\Delta\phi+\gamma\lambda(\theta-\phi).
\end{aligned}
\end{align}
where stream function $\psi$ satisfies
\begin{align}\label{b3}
v_{1} = -\frac{\partial \psi}{\partial z},\ \ \ \
v_{2} = \frac{\partial \psi}{\partial x},
\end{align}
and the nonlinear part reads
\begin{align}\label{b3,}
J\left(\psi,\theta\right)
= \frac{\partial \psi}{\partial x}\frac{\partial \theta}{\partial z}
 - \frac{\partial \psi}{\partial z}\frac{\partial \theta}{\partial x}.
\end{align}

The non-dimensionalization procedure gives rise to the following seven numbers: the Darcy-Rayleigh number $Ra = \frac{\rho_{0} g \beta \left( T_{l} - T_{u} \right)  \left(\rho c\right)_{f} K }{ \epsilon \mu_{f} \kappa_{f}}$, the dimensionless heat transfer coefficient $\lambda = \frac{h }{\epsilon \kappa_{f}}$, the modified conductivity $\gamma = \frac{\epsilon \kappa_{f}}{ \left( 1 - \epsilon \right)\kappa_{s}}$, the diffusion ration $\alpha = \frac{\left(\rho c\right)_{s}}{\left(\rho c\right)_{f} } \frac{\kappa_{f}}{\kappa_{s}}$, the Darcy number $Da = K$, the Prandtl number $Pr = \frac{ \mu_{f} \epsilon \left(\rho c\right)_{f}}{\rho_{0} \kappa_{f}}$ and the couple stress number $C = \frac{\mu_{c}}{\mu_{f}}$.
%\begin{align}\label{b5}
%\begin{aligned}
%& Ra = \frac{\rho_{0} g \beta \left( T_{l} - T_{u} \right)  \left(\rho c\right)_{f} K }{ \epsilon \mu_{f} \kappa_{f}},
%\quad
%\gamma = \frac{\epsilon \kappa_{f}}{ \left( 1 - \epsilon \right)\kappa_{s}},
%\quad
%\lambda = \frac{h }{\epsilon \kappa_{f}},
%\quad
%\alpha = \frac{\left(\rho c\right)_{s}}{\left(\rho c\right)_{f} } \frac{\kappa_{f}}{\kappa_{s}},
%\\
% &C = \frac{\mu_{c}}{\mu_{f} },
%\quad
%Pr = \frac{ \mu_{f} \epsilon \left(\rho c\right)_{f}}{\rho_{0} \kappa_{f} },
%\quad
%Da = K,
%\quad
%J\left(\psi,\theta\right) = \frac{\partial \psi}{\partial x}\frac{\partial \theta}{\partial z}
%- \frac{\partial \psi}{\partial z}\frac{\partial \theta}{\partial x}.
%\end{aligned}
%\end{align}
Now, for the bounded set $\Omega\subset R^2$, we give the boundary and initial conditions to the equation \eqref{b4} as follows:
\begin{align}\label{b6}
\begin{aligned}
\psi|_{\partial \Omega} &=\Delta\psi|_{\partial \Omega}=0,
\\
\theta|_{\partial \Omega} &=\phi|_{\partial \Omega}=0,
\\
\mathbf{V}(0,x,z) &= \left(\psi(0,x,z),\theta(0,x,z),\phi(0,x,z) \right).
\end{aligned}
\end{align}

\section{Well-posedness  of solution}

\subsection{Mathematical setting}
We denote $\langle\cdot,\cdot\rangle$ as the inner product of $L^2$ and $\|\cdot\|$ as $L^2$-norm. Let $\mathbf{V}=(\psi,\theta,\phi)$ and introduce the following spaces
\begin{align}
\begin{aligned}
 \mathbf{X}
&=\{\mathbf{V}\in C^{\infty}(\Omega)\times C^{\infty}(\Omega)\times C^{\infty}(\Omega) |~\mathbf{V} ~\textrm{satifies}~ \eqref{b6}\},
\\
 \mathbf{Y}
&=H^{2}(\Omega)\times L^2(\Omega)\times L^2(\Omega),
\\
 \mathbf{Y_{\frac{1}{2}}}
&=\{\mathbf{V}\in H^3(\Omega)\times H^1(\Omega)\times H^1(\Omega)
| ~\mathbf{V} ~\textrm{satifies}~ \eqref{b6}\},
\\
 \mathbf{Y_{1}}
&=\{\mathbf{V}\in H^4(\Omega)\times H^2(\Omega)\times H^2(\Omega)|~\mathbf{V} ~\textrm{satifies}~ \eqref{b6}\},
\\
 \mathbf{Y_{\frac{k}{2}}}
&=\{\mathbf{V}\in H^{k+2}(\Omega)\times H^{k}(\Omega)\times H^{k}(\Omega)|~\mathbf{V} ~\textrm{satifies}~ \eqref{b6}\}
\ \ (k\geq 2).
\end{aligned}
\end{align}

Because $\mathbf{X}\subset \mathbf{Y_{\frac{1}{2}}}$ is dense and $\mathbf{Y_{\frac{1}{2}}}$ is separable, there exists a set of orthonormal basis $\{\tilde{e}_{m}\}\times\{e_{m}\}\times\{e_{m}\}$ which satisfies the following conditions
\begin{align}\label{w1}
\begin{aligned}
& ~\Delta \tilde{e}_{m}=\tilde{\mu}_{m}\tilde{e}_{m},\ \ ~ \ \  \Delta e_{m}=\mu_{m} e_{m},
\\
& \langle\tilde{e}_{m_{1}},\tilde{e}_{m_{2}} \rangle=\langle e_{m_{1}},e_{m_{2}}\rangle=\delta_{m_{1}m_{2}}.
\end{aligned}
\end{align}

In addition, according to \cite{Grisvard2011}, when $\psi\in H^{2}\cap H_{0}^{1}$, there exists a constant $c>0$(depends on $\Omega$) such that 
\begin{align}\label{constant1}
  \|\psi\|_{H^{2}}\leq c\|\Delta\psi\|,
\end{align}
then by \textbf{Poincar\'e inequality}, $\|\nabla\Delta\psi\|$ is equivalent to $\|\psi\|_{H^{3}}$ when $\mathbf{V}\in \mathbf{Y_{\frac{1}{2}}}$.
\subsection{Existence and uniqueness of the global solution}
Before giving the global solution of the equation \eqref{b4}, we consider the existence and uniqueness of weak solution.
\begin{definition}(\textbf{Weak solution})
We say $\mathbf{V}=(\psi,\theta,\phi)\in L^{\infty}((0,T),\mathbf{Y})\cap L^2((0,T),\mathbf{Y_{\frac{1}{2}}})$ is the \textbf{weak solution} of the equation \eqref{b4}, if for any $ \mathbf{W}=(w_{1},w_{2},w_{3})\in X$ and $0<t<T$ $(T>0$ is arbitrary fixed constant), we have
\begin{align}\label{w2}
\begin{aligned}
 & \frac{Da}{Pr}\langle\Delta\psi(t),w_{1}\rangle + \langle\theta(t),w_{2}\rangle + \alpha\langle\phi(t),w_{3}\rangle
\\
= &\int_{0}^{t}[-C\langle\nabla\Delta\psi,\nabla w_{1}\rangle-\left\langle \Delta\psi, w_{1}\right\rangle  + Ra\langle\frac{\partial\theta}{\partial x},w_{1}\rangle
    -\langle\nabla\theta,\nabla w_{2}\rangle 
\\
 &+\lambda\langle\phi-\theta,w_{2}\rangle -\langle J(\psi,\theta),w_{2} \rangle -\langle\nabla\phi,\nabla w_{3}\rangle
   +\gamma\lambda\langle\theta -\phi,w_{3}\rangle ] d\tau
\\
 & +\frac{Da}{Pr}\langle\Delta\psi(0),w_{1}\rangle + \langle\theta(0),w_{2}\rangle +\alpha\langle\phi(0),w_{3}\rangle.
\end{aligned}
\end{align}
\end{definition}
According to the Galerkin method, we can get the following results.

\begin{theorem}\label{w3}(\textbf{Existence of weak solution}) For any initial value $\mathbf{V}(0)\in \mathbf{Y_{\frac{1}{2}}}$, the equation \eqref{b4} has a weak solution $\mathbf{V}(t)\in L_{loc}^{\infty}((0,\infty),\mathbf{Y})\cap L_{loc}^{2}((0,\infty),\mathbf{Y_{\frac{1}{2}}}).$
\end{theorem}
\begin{proof}\textbf{This proof is divided into three steps. }

\textbf{Step 1: Construct an approximate solution.}
Let $\mathbf{V^{m}}=(\psi^m,\theta^m,\phi^m)(m\in \mathbb{N})$, where
\begin{align*}
 \psi^{m} &= \sum\limits_{i=1 }^{m} \psi_{i }^{m} \left(t\right) \tilde{e}_{i },
\\
 \theta^{m} &= \sum\limits_{i=1 }^{m} \theta_{i }^{m} \left(t\right) e_{i },
\\
 \phi^{m} &= \sum\limits_{i=1 }^{m} \phi_{i }^{m} \left(t\right) e_{i }.
\end{align*}
%and $\psi_{i}^{m}(t), \theta_{i}^{m}(t)$,$\phi_{i}^{m}(t)\in C^{\infty}(\Omega).$
Setting $\mathbf{V}=\mathbf{V^m}$ in the equation \eqref{b4} and taking the inner product of equation \eqref{b4} with $(\tilde{e_{i}},e_{i},e_{i})(i=1,\cdots,m)$, we have the set of the following ordinary equations
\begin{align}\label{w4}
\begin{aligned}
 \frac{Da}{Pr}\frac{d\psi_{i}^{m}(t)}{dt}
& =(C\tilde{\mu_{i}}-1)\psi_{i}^{m}(t)+Raf_{1i}(t),
\\
 \frac{d\theta_{i}^{m}(t)}{dt}
& =(\mu_{i}-\lambda)\theta_{i}^{m}(t)+\lambda\phi_{i}^{m}(t)-f_{2i}(t),
\\
 \alpha\frac{ d\phi_{i}^{m}(t)}{dt}
 & =(\mu_{i}-\gamma\lambda)\phi_{i}^{m}(t)+\gamma\lambda\theta_{i}^{m}(t),
\end{aligned}
\end{align}
where
\begin{align*}
f_{1i}& =\sum\limits_{k=1}^{m}\frac{\theta_{k}^{m}(t)}{\tilde{\mu}_{k}}\langle\frac{\partial e_{k}}{\partial x}, \tilde{e_{i}}\rangle,
\\
f_{2i}& =\sum\limits_{k_{1},k_{2}=1}^{m}\psi_{k_{1}}^{m}(t)\theta_{k_{2}}^{m}(t)\langle\frac{\tilde{\partial e_{k_{1}}}}{\partial x}\frac{\partial e_{k_{2}}}{\partial z}-\frac{\tilde{\partial e_{k_{1}}}}{\partial z}\frac{\partial e_{k_{2}}}{\partial x},e_{i}\rangle.
\end{align*}

Based on the fundamental theory of ordinary differential equations in the space of smooth functions, the equation \eqref{w4} has a smooth local solution , by which we have
\begin{align}\label{w5}
\begin{aligned}
  &\frac{Da}{Pr}\langle\Delta\psi^{m},w_{1}\rangle + \langle\theta^{m},w_{2}\rangle + \alpha\langle\phi^{m},w_{3}\rangle
\\
= & \int_{0}^{t} [-C\langle\nabla\Delta\psi^{m},\Delta w_{1}\rangle-\left\langle \Delta\psi^{m},w_{1}\right\rangle  + Ra\langle\frac{\partial\theta^{m}}{\partial x},w_{1}\rangle
     -\langle\nabla\theta^{m},\nabla w_{2}\rangle
\\
 & +\lambda\langle\phi^{m}-\theta^{m},w_{2}\rangle-\langle J(\psi^{m},\theta^{m}),w_{2}\rangle
  -\langle\nabla\phi^{m},\nabla w_{3}\rangle+\gamma\lambda\langle\theta^{m}-\phi^{m},w_{3}\rangle ]d\tau
\\
 & +\frac{Da}{Pr}\langle\Delta\psi^{m}(0),w_{1}\rangle  +\langle\theta^{m}(0),w_{2}\rangle
   +\alpha\langle\phi^{m}(0),w_{3}\rangle
\end{aligned}
\end{align}
for any $ 0<t<T$.

\textbf{Step 2: Priori estimates.}
Substituting $\mathbf{V^{m}}=(\Delta\psi^m,\theta^m,\phi^m)(m\in \mathbb{N})$ into the equation \eqref{b4} and taking the inner product of with $\mathbf{V^{m}}$, we can obtain
\begin{align}\label{w5,}
\begin{aligned}
 & \frac{1}{2} \frac{d}{dt} \left( \frac{Da}{Pr}\left\|\Delta \psi^{m} \right\|^{2} + \left\|\theta^{m} \right\|^{2}
        + \alpha\left\|\phi^{m} \right\|^{2} \right)
\\
=& - \left(C \left\|\nabla\Delta \psi^{m} \right\|^{2}
   + \left\|\nabla\theta^{m} \right\|^{2}
   +  \left\|\nabla\phi^{m} \right\|^{2} \right)
\\
& - \left\|\Delta \psi^{m} \right\|^{2} -Ra\left\langle \theta,\frac{\partial\Delta\psi^{m}}{\partial x}\right\rangle 
\\
&  +\lambda\left\langle \phi^{m}-\theta^{m},\theta^{m}\right\rangle
+\gamma\lambda\left\langle \theta^{m}-\phi^{m},\phi^{m}\right\rangle.  
\end{aligned}
\end{align}
Then, the Young inequality is used to make estimates
\begin{align}\label{w5,1}
 |Ra\langle\theta^{m}, \frac{\partial \Delta\psi^{m}}{\partial x}\rangle|
\leq\frac{C}{2} \|\nabla\Delta\psi^{m} \|^2 +\frac{Ra^2}{2C}\|\theta^{m}\|^2,
\end{align}
and
\begin{align}\label{w5,2}
  | (\lambda + \gamma\lambda) \langle\phi^{m},\theta^{m}\rangle|
\leq (\gamma\lambda + \frac{\lambda}{4})\|\phi^{m}\|^2
     +(\frac{\gamma\lambda}{4} + \lambda)\|\theta^{m}\|^2,
\end{align}

In combination with \eqref{w5,}-\eqref{w5,2}, we have
\begin{align}\label{b2.8}
 \frac{d}{dt} \left( \frac{Da}{Pr}\left\|\Delta \psi^{m} \right\|^{2} + \left\|\theta^{m} \right\|^{2}
             + \alpha\left\|\phi^{m} \right\|^{2} \right)
\leq -M_{1}\|\mathbf{V}^{m}\|_{\mathbf{Y_{\frac{1}{2}}}}^{2}+M_{2}\|\mathbf{V}^{m}\|_{\mathbf{Y}},
\end{align}
where
\begin{align*}
M_{1} =2\min\left\{\frac{PrC}{2Da},1,\frac{1}{\alpha}\right\},
\ \ \ \
M_{2} =2\max\left\{1,\frac{Ra^2}{2C}+\frac{\lambda\gamma}{4},\frac{\lambda}{4\alpha}\right\},
\end{align*}
According to \textbf{Gronwall inequality}, it yields that
\begin{align*}
  \|\mathbf{V^{m}}\|_{\mathbf{Y}}^2
\leq \|\mathbf{V^{m}}(0)\|_{\mathbf{Y}}^2e^{M_{2}t},
\end{align*}
and
\begin{align*}
  \int_{0}^{t} \|\mathbf{V^{m}}\|_{\mathbf{Y_{\frac{1}{2}}}}^2 d\tau
\leq \frac{e^{M_{2}t}}{M_{1}}\|\mathbf{V^{m}}(0)\|_{\mathbf{Y}}^2 ,
\end{align*}
which shows the approximate solution
\[
\mathbf{V^{m}}\in L^{\infty}((0,T),\mathbf{Y})
\cap L^2((0,T),\mathbf{Y_{\frac{1}{2}}})\quad
\text{for} \quad \forall ~m\in \mathbb{N}.
\]

\textbf{Step 3: Getting a convergent sequence in some sense.}
In order to take the limit of the equation \eqref{w5}, the approximate solution must be convergent in some sense. Firstly, we prove $\{\psi_{i}^{m}(t)\}$, $\{\theta_{i}^{m}(t)\}$ and $\{\phi_{i}^{m}(t)\}$ are equicontinuous for fixed $i$ and $m=1,2,\cdots.$

On one hand, we have
\begin{align*}
|\psi_{i}^{m}(t+h)-\psi_{i}^{m}(t)|
&= |\int_{t}^{t+h}\frac{d\psi_{i}^{m}(\tau)}{d\tau}d\tau |
\\
&= |\int_{t}^{t+h} \left[\langle\frac{Pr}{Da}(C\Delta-1)\Delta\psi^{m},\tilde{e_{i}}\rangle
  +\frac{PrRa}{Da}\langle\frac{\partial \theta^m}{\partial x},\tilde{e_{i}}\rangle\right] d\tau |
  \\
&\leq M_{3}h^{\frac{1}{2}},
\end{align*}
where
\begin{align*}
\begin{aligned}
M_{3}=3\max_{0<t<T} \{
&\frac{PrC}{Da}\|\nabla \tilde{e_{i}}\|[\int_{t}^{t+h}\|\nabla\Delta\psi^m\|^2 d\tau]^{\frac{1}{2}},~ \frac{Pr \|\tilde{e_{i}}\|}{Da}\|\Delta\psi^{m}\|,~
\\
& 
\frac{PrRa}{Da}\|\nabla\tilde{e_{i}}\|[\int_{t}^{t+h}\|\theta^m\|^2 d\tau]^{\frac{1}{2}}
\}.
\end{aligned}
\end{align*}
On other hand, we have
\begin{align*}
|\theta_{i}^{m}(t+h)-\theta_{i}^{m}(t)|\leq M_{4}h^{\frac{1}{2}},
\end{align*}
and
\begin{align*}
|\phi_{i}^{m}(t+h)-\phi_{i}^{m}(t)|\leq M_{5}h^{\frac{1}{2}},
\end{align*}
where
\begin{align*}
M_{4}=3\max_{0<t<T}\{
  &
  \|\nabla e_{i}\|[\int_{t}^{t+h}\|\nabla\theta^m\|^2 d\tau]^{\frac{1}{2}},~ \lambda[\|\phi^m\|+\|\theta^m\|]\|e_{i}\|,
  \\
  &2cM_{so}\|\Delta\psi^m\|\|e_{i}\|[\int_{t}^{t+h}\|\nabla\theta^m\|^2 d\tau]^{\frac{1}{2}}
\},
\end{align*}
\begin{align*}
M_{5}=2\max_{0<t<T}\{\frac{1}{\alpha}\|\nabla e_{i}\|[\int_{t}^{t+h}\|\nabla\phi^{m}\|^{2}d\tau]^{\frac{1}{2}},~\frac{\gamma\lambda}{\alpha}[\|\theta^m\|+\|\phi^m\|]\|e_{i}\|\},
\end{align*} 
where $c$ is in \eqref{constant1} and $M_{so}$ is from \textbf{Sobolev inequality}. 

Note that $\{\psi_{i}^{m}(t)\}$, $\{\theta_{i}^{m}(t)\}$ and $\{\phi_{i}^{m}(t)\}$ are uniformly bounded by the results given in \textbf{Step 2}.
Following the \textbf{Ascoli-Arzela} theorem, there exists convergent subsequences of $\{\psi^{m}_{i}(t)\}$, $\{\theta^{m}_{i}(t)\}$ and $\{\phi^{m}_{i}(t)\}$, denote by $\{\psi^{m}_{i}(t)\}$, $\{\theta^{m}_{i}(t)\}$ and $\{\phi^{m}_{i}(t)\}$, such that
$$
\psi_{i}^{m}\rightarrow \psi_{i}, \ \ \theta_{i}^{m}(t)\rightarrow\theta_{i}(t),\ \ \phi_{i}^{m}(t)\rightarrow \phi_{i}(t)
$$
as $m\rightarrow\infty$ uniformly, where $\psi_{i}(t)$, $\theta_{i}(t)$ and $\phi_{i}(t)$ are all continuous.

Denote $$
\mathbf{V_{0}}
=\left(\sum\limits_{i=1}^{+\infty}\psi_{i}(t)\tilde{e_{i}},
  \sum\limits_{i=1}^{+\infty}\theta_{i}(t)e_{i},
  \sum\limits_{i=1}^{+\infty}\phi_{i}(t)e_{i}\right),
$$
then for any $\tilde{\mathbf{V}}\in X$ we have
\begin{align}\label{w6}
\lim\limits_{m\rightarrow\infty}\sup\limits_{0\leq t\leq T}\langle \mathbf{V^{m}}-\mathbf{V_{0}},\tilde{\mathbf{V}} \rangle=0,
\end{align}
and infer that
\begin{align}\label{w7}
\lim\limits_{m\rightarrow\infty}\int_{0}^{T}|\langle\mathbf{V^m}-\mathbf{V_{0}},\mathbf{\tilde{V}} \rangle|^2d\tau=0.
\end{align}

Thus, it means $\{\mathbf{V^m}\}$ is uniformly weakly convergent to $\mathbf{V_{0}}$ in $L^2((0,T),\mathbf{Y_{\frac{1}{2}}})$. From (lemma 4.28 in \cite{2011}), we get
\begin{align}
\begin{cases}
D^{\alpha}\psi^{m}\rightarrow D^{\alpha}\psi_{0}, \ \  \ \  \forall ~|\alpha|\leq 2,
\\
\theta^{m}\rightarrow\theta_{0},
\\
\phi^m\rightarrow\phi_{0},
\end{cases}
\end{align}
in $L^2((0,T)\times\Omega)$ as $m\rightarrow\infty$.

For simplicity, we only focus on the nonlinear term. A simple calculation shows that
\begin{align*}
 &|\int_{0}^{t}\langle\frac{\partial \psi^m}{\partial x}\frac{\partial \theta^m}{\partial z}
 -\frac{\partial \psi^m}{\partial z}\frac{\partial \theta^m}{\partial x}-(\frac{\partial\phi_{0}}{\partial x}\frac{\partial \theta_{0}}{\partial z}-\frac{\partial \phi_{0}}{\partial z}\frac{\partial\theta_{0}}{\partial x}),w_{2}\rangle d\tau|
\\
\leq& |\int_{0}^{t}\langle(\frac{\partial^2 \psi^m}{\partial x\partial z}-\frac{\partial^2\psi_{0}}{\partial x \partial z})\theta^m,w_{2}\rangle +\left\langle (\frac{\partial\psi^m}{\partial x}-\frac{\partial\psi_{0}}{\partial x})\theta^m,\frac{\partial w_{2}}{\partial z}\right\rangle  d\tau|
\\
& +|\int_{0}^{t}\langle\frac{\partial^2 \psi_{0}}{\partial x\partial z}(\theta^m-\theta_{0})),w_{2}\rangle +\left\langle \frac{\partial\psi_{0}}{\partial x}(\theta^m-\theta_{0}),\frac{\partial w_{2}}{\partial z}\right\rangle d\tau|
\\
&|\int_{0}^{t}\langle(\frac{\partial^2 \psi_{0}}{\partial x\partial z}-\frac{\partial^2\psi^{m}}{\partial x \partial z})\theta_{0},w_{2}\rangle +\left\langle (\frac{\partial\psi_{0}}{\partial z}-\frac{\partial\psi^{m}}{\partial z})\theta_{0},\frac{\partial w_{2}}{\partial x}\right\rangle  d\tau|
\\
&+|\int_{0}^{t}\langle\frac{\partial^2 \psi^{m}}{\partial x\partial z}(\theta_{0}-\theta^{m})),w_{2}\rangle +\left\langle \frac{\partial\psi^{m}}{\partial z}(\theta_{0}-\theta^{m}),\frac{\partial w_{2}}{\partial x}\right\rangle d\tau|
\\
&
\longrightarrow 0~(m\rightarrow \infty).
\end{align*}

Taking the limit on both sides of the equation \eqref{w5}, we obtain
\begin{align*}
\begin{aligned}
  & \frac{Da}{Pr}\langle\Delta\psi_{0},w_{1}\rangle + \langle\theta_{0},w_{2}\rangle + \alpha\langle\phi_{0},w_{3}\rangle
\\
= & \int_{0}^{t}  [-C\langle( \nabla\Delta \psi_{0},\nabla w_{1}\rangle-\left\langle \Delta\psi_{0},w_{1}\right\rangle 
    + Ra\langle\frac{\partial\theta_{0}}{\partial x},w_{1}\rangle -\langle\nabla\theta_{0},\nabla w_{2}\rangle
\\
  &  + \lambda\langle\phi_{0} - \theta_{0},w_{2}\rangle
    -\langle J(\psi_{0},\theta_{0}),w_{2}\rangle -\langle\nabla\phi_{0},\nabla w_{3}\rangle
    + \gamma\lambda\langle\theta_{0}-\phi_{0},w_{3}\rangle ] d\tau
\\
  & +\frac{Da}{Pr}\langle\Delta\psi_{0}(0),w_{1}\rangle + \langle\theta_{0}(0),w_{2}\rangle
    +\alpha\langle\phi_{0}(0),w_{3}\rangle,
\end{aligned}
\end{align*}
which means $\mathbf{V_{0}}$ is the weak solution of the equation \eqref{b4}.
\end{proof}

Moreover, by the two-dimensional Sobolev embedding theorem\cite{2011} (page 21), we have the following theorem.

\begin{theorem}\label{w8}(\textbf{Uniqueness}) The weak solution of the equation \eqref{b4} is unique.
\end{theorem}
\begin{proof}
Assume both $\mathbf{V_{1}}=(\psi_{1}, \theta_{1}, \phi_{1})$ and $\mathbf{V_{2}}=(\psi_{2}, \theta_{2}, \phi_{2})$ are weak solutions of the equation \eqref{b4}.
Let $\tilde{\mathbf{V}}= (\tilde{\psi}, \tilde{\theta}, \tilde{\phi}) =\mathbf{V_{1}}-\mathbf{V_{2}}$, that is
$$
\tilde{\psi}=\psi_{1}-\psi_{2}, \ \
\tilde{\theta}=\theta_{1}-\theta_{2},\ \
\tilde{\phi}=\phi_{1}-\phi_{2}.
$$
By simple calculation, we know $\mathbf{\tilde{V}(0)}=0$ and
$$
J(\psi_{1},\theta_{1})-J(\psi_{2},\theta_{2})=J(\tilde{\psi},\theta_{1})+J(\psi_{2},\tilde{\theta}),
\ \ \ \
\langle J(\psi_{2},\tilde{\theta}),\tilde{\theta}\rangle=0.
$$

Putting $ \mathbf{W}=\mathbf{\tilde{V}}$ into the equation \eqref{w2}, we obtain that
\begin{align*}
 &\frac{Da}{Pr}\|\nabla\tilde{\psi}\|^2+\|\tilde{\theta}\|^2+\alpha\|\tilde{\phi}\|^2
\\
= &\int_{0}^{t}[-(C\|\Delta\tilde{\psi}\|^2+\|\nabla\tilde{\theta}\|^2+\|\nabla\tilde{\phi}\|^2)
    -(\|\nabla\tilde{\psi}\|^2+\lambda\|\tilde{\theta}\|^2 +\gamma\lambda\|\tilde{\phi}\|^2)
\\
& - Ra\langle\tilde{\theta},\frac{\partial\tilde{\psi}}{\partial x}\rangle
  +(\lambda + \gamma\lambda)\langle\tilde{\phi},\tilde{\theta}\rangle
  -\langle J(\tilde{\psi},\theta_{1}),\tilde{\theta}\rangle ]d\tau.
\end{align*}
By \textbf{Young inequalities}, we receive the following estimation
\begin{align*}
  |Ra\langle \tilde{\theta},\frac{\partial\tilde{\psi}}{\partial x}\rangle|
\leq \|\nabla\tilde{\psi}\|^2 + \frac{Ra^2}{4}\|\tilde{\theta}\|^2,
\end{align*}
and
\begin{align*}
  |(\lambda+\gamma\lambda)\langle\tilde{\phi},\tilde{\theta}\rangle|
\leq (\gamma\lambda+\frac{\lambda}{4})\|\tilde{\phi}\|^2
     +(\frac{\gamma\lambda}{4}+\lambda)\|\tilde{\theta}\|^2,
\end{align*}
Besides, on the basis of \textbf{Solobev embedding theorem}, we have
\begin{align*}
 2\int_{\Omega}|\nabla\tilde{\psi}\|\nabla\theta_{1}\|\tilde{\theta}|dxdz
& \leq  2\|\tilde{\theta}\|_{L^{\infty}} | \int_{\Omega}|\nabla\tilde{\psi}||\nabla\theta_{1}|dxdz |
\\
& \leq \|\nabla\tilde{\theta}\|^2 + M_{so}^2\|\nabla\theta_{1}\|^2\|\nabla\tilde{\psi}\|^2,
\end{align*}
where $M_{so}$ just depends on $\Omega$.

To sum up, we come to the conclusion
\begin{align*}
\frac{Da}{Pr}\|\nabla\tilde{\psi}\|^2+\|\tilde{\theta}\|^2+\alpha\|\tilde{\phi}\|^2\leq \int_{0}^{t}\alpha(\tau)(\frac{Da}{Pr}\|\nabla\tilde{\psi}\|^2+\|\tilde{\theta}\|^2+\alpha\|\tilde{\phi}\|^2)d\tau,
\end{align*}
where nonnegative function $\alpha(\tau)=\max\{\frac{M_{so}^2\|\nabla \theta_{1}\|^2 Pr}{Da},\frac{Ra^2+\gamma\lambda}{4},\frac{\lambda}{4\alpha}\}\in L^{1}(0,T)$.
Finally, according to \textbf{Gronwall inequality}, we get
$$
\frac{Da}{Pr}\|\nabla\tilde{\psi}\|^2+\|\tilde{\theta}\|^2+\alpha\|\tilde{\phi}\|^2\leq0,
$$
which means the weak solution is unique.

\end{proof}

Next, we will improve the regularity of weak solution with respect to $t$. Let us start with the definition of the global solution.
\begin{definition}(\textbf{Global solution}) We say $\mathbf{V}$ is the global solution of the equation\eqref{b4}, if for any $\mathbf{W}=(w_{1},w_{2},w_{3})\in \mathbf{X}$, there is
\begin{align}\label{g1}
\begin{aligned}
 & \frac{Da}{Pr}\langle\frac{d\Delta\psi(t)}{dt},w_{1}\rangle + \langle\frac{d\theta(t)}{dt},w_{2}\rangle
 + \alpha\langle\frac{d\phi(t)}{dt},w_{3}\rangle
\\
= & -C\langle\nabla\Delta\psi,\nabla w_{1}\rangle-\left\langle \Delta\psi,w_{1}\right\rangle  + Ra\langle\frac{\partial\theta}{\partial x},w_{1}\rangle
    - \langle\nabla\theta,\nabla w_{2}\rangle + \lambda\langle\phi - \theta,w_{2}\rangle
\\
 & -\langle J(\psi,\theta),w_{2}\rangle - \langle\nabla\phi,\nabla w_{3}\rangle
   + \gamma\lambda\langle\theta-\phi,w_{3}\rangle.
\end{aligned}
\end{align}
\end{definition}

For the equation \eqref{b4}, we proved the following result.
\begin{theorem}(\textbf{Existence and uniqueness of global solution}) If the initial value $\mathbf{V(0)}\in \mathbf{Y_{1}}$, the equation \eqref{b4} has a unique global solution
$$
\mathbf{V}\in W_{loc}^{1,\infty}((0,\infty),\mathbf{Y})\cap W_{loc}^{1,2}((0,\infty),\mathbf{Y_{\frac{1}{2}}}).
$$
\end{theorem}

\begin{proof}\textbf{Firstly, we will estimate the approximate solution $\mathbf{V^{m}}$ which we have got.}

Obviously, $\mathbf{V^{m}}=(\psi^m,\theta^m,\phi^m)$ satisfies the following equations
\begin{align}\label{g2}
\begin{aligned}
 &\int_{0}^{t} ( \langle \Delta\ddot{\psi}^m,w_{1}\rangle + \langle\ddot{\theta}^{m},w_{2}\rangle
   +  \langle \ddot{\phi}^{m},w_{3}\rangle ) d\tau
\\
=& \int_{0}^{t} [ -\frac{PrC}{Da}\langle\nabla\Delta\dot{\psi}^m,\nabla w_{1} \rangle-\frac{Pr}{Da}\left\langle \Delta\dot{\psi}^m,w_{1}\right\rangle 
     + \frac{Pr Ra}{Da} \langle\frac{\partial \dot{\theta}^{m}}{\partial x},w_{1}\rangle
     - \langle\nabla\dot{\theta}^{m},\nabla w_{2}\rangle
\\
 &
   + \lambda\langle \dot{\phi}^{m}-\dot{\theta}^{m} , w_{2} \rangle - \langle J(\dot{\psi}^{m},\theta^{m}) + J( \psi^{m},\dot{\theta}^{m}) ,w_{2} \rangle
    - \frac{1}{\alpha}\langle\nabla\dot{\phi}^{m},\nabla w_{3}\rangle 
    + \frac{\gamma\lambda}{\alpha}\langle\dot{\theta}^{m} - \dot{\phi}^{m}, w_{3} \rangle] d\tau.
\end{aligned}
\end{align}
Let $(w_{1},w_{2},w_{3})=(\Delta\dot{\psi}^{m},\dot{\theta}^{m},\dot{\phi}^m)$, then we have
\begin{align}\label{g3}
\begin{aligned}
  & \|\Delta\dot{\psi}^{m}\|^2 + \|\dot{\theta}^{m}\|^2 + \|\dot{\phi}^{m}\|^2
\\
= & -2 \int_{0}^{t} [\frac{PrC}{Da}\|\nabla\Delta\dot{\psi}^{m}\|^2+\frac{Pr}{Da}\|\Delta\dot{\psi}^{m}\|^2
    - \frac{ Pr Ra}{Da} \langle\frac{\partial\dot{\theta}^{m}}{\partial x},\Delta\dot{\psi}^{m} \rangle
    + \|\nabla\dot{\theta}^m\|^2
\\
 & + \lambda\langle\dot{\theta}^{m} - \dot{\phi}^{m},\dot{\theta}^{m} \rangle
   + \langle J(\dot{\psi}^{m},\theta^{m}),\dot{\theta}^{m} \rangle
   + \frac{1}{\alpha}\|\nabla\dot{\phi}^{m}\|^2 
   + \frac{\gamma\lambda}{\alpha} \langle\dot{\phi}^{m} -\dot{\theta}^{m},\dot{\phi}^{m} \rangle ] d\tau
\\
&  + \|\Delta\dot{\psi}^{m}(0)\|^2 + \|\dot{\theta}^{m}(0)\|^2 +\|\dot{\phi}^{m}(0)\|^2.
\end{aligned}
\end{align}
Because 
\begin{align}\label{g4}
\begin{aligned}
  \frac{Da}{Pr}\langle\Delta\dot{\psi}^{m},w_{1}\rangle
= & \langle(C\Delta-1)\Delta\psi^{m}, w_{1}\rangle+Ra\langle\frac{\partial\theta^{m}}{\partial x}, w_{1} \rangle,
\\
  \langle\dot{\theta}^{m},w_{2}\rangle
= & \langle \Delta\theta^{m},w_{2}\rangle + \lambda\langle\phi^m - \theta^m,w_{2}\rangle
  -\langle J(\psi^{m},\theta^{m}),w_{2}\rangle,
\\
  \alpha\langle \dot{\phi}^{m},w_{3}\rangle
= &\langle\Delta\phi^{m},w_{3}\rangle + \gamma\lambda\langle\theta^{m}-\phi^m,w_{3}\rangle,
\end{aligned}
\end{align}
let $(w_{1},w_{2},w_{3})=(\Delta\dot{\psi}^{m},\dot{\theta}^{m},\dot{\phi}^{m})$ in \eqref{g4}, and by \textbf{Young inequality}, we get
\begin{align}\label{g5}
  \|\Delta\dot{\psi}^{m}(0)\|^{2}+\|\dot{\theta}^{m}(0)\|^2+\|\dot{\phi}^{m}(0)\|^{2}
\leq g \left(\psi^{m}(0),\theta^{m}(0),\phi^{m}(0)\right),
\end{align}
where
\begin{align*}
  g(\psi^{m}(0),\theta^{m}(0),\phi^{m}(0))
=& \frac{3Pr^{2}C^2}{Da^2}\|\Delta^2\psi^{m}(0)\|^2+\frac{3Pr^2}{Da^2}\|\Delta\psi^{m}(0)\|^2+3\|\Delta\theta^{m}(0)\|^2
\\
& +\frac{2}{\alpha^2}\|\Delta\phi^{m}(0)\|^2+ 12\|\nabla\psi^{m}(0)\nabla\theta^{m}(0)\|^2+\frac{3Pr^2 Ra^2}{Da^2}\|\nabla\theta^m(0)\|^2
\\
&  +(12\lambda^2+\frac{8\gamma^2\lambda^2}{\alpha^2})(\|\theta^{m}(0)\|^2+\|\phi^m(0)\|^2).
\end{align*}
%Because $\mathbf{Y_{\frac{3}{2}}}\subset\mathbf{Y_{1}}$ is dense and compact, when $\mathbf{V}(0)\in\mathbf{Y_{\frac{3}{2}}}$ and $\mathbf{V}^{m}(0)\rightarrow\mathbf{V}(0)$ in $\mathbf{Y_{\frac{3}{2}}}$, we have $\mathbf{V}^{m}(0)\rightarrow \mathbf{V}(0)$ in $\mathbf{Y_{1}}$. 

Thus, we choose $\mathbf{V}(0)\in\mathbf{Y_{1}}$, then $\{\|\mathbf{V}^{m}(0)\|_{\mathbf{Y_{1}}}\}$ is bounded and $g(\psi^{m}(0),\theta^{m}(0),\phi^{m}(0))$ is bounded.

\textbf{Secondly, we can obtain the solution by functional analysis.}
According to the \textbf{Sobolev inequality} and \eqref{constant1}, we know
\begin{align*}
  |-2\int_{0}^{t}\langle J(\dot{\psi}^{m},\theta^{m}),\dot{\theta}^{m}\rangle d\tau|
& \leq 4\int_{0}^{t}\|\nabla\theta^{m}\|\|\nabla\dot{\psi}^{m}\dot{\theta}^{m}\|d\tau
\\
& \leq 4cM_{so}\int_{0}^{t}\|\nabla\theta^{m}\|\|\dot{\mathbf{V}}^m\|^2_{\mathbf{Y}}d\tau,
\end{align*}
where $M_{so}>0$ is from \textbf{Sobolev inequality} and $c>0$ is in \eqref{constant1}. And by \textbf{Young inequality}, we receive that
\begin{align*}
  |-2\frac{PrRa}{Da}\int_{0}^{t}\langle \frac{\partial \dot{\theta}^{m}}{\partial x},\Delta\dot{\psi}^{m}\rangle d\tau|
\leq \frac{PrRa^2}{3Da}\int_{0}^{t}\|\dot{\theta}^{m}\|^2d\tau + \frac{3Pr}{Da}\int_{0}^{t}\|\Delta\dot{\psi}^{m}\|^2 d\tau,
\end{align*}
and
\begin{align*}
|-2\lambda\int_{0}^{t}\langle\dot{\theta}^{m}-\dot{\phi}^{m},\dot{\theta}^{m}\rangle
  -\frac{2\gamma\lambda}{\alpha}\langle\dot{\phi}^{m}-\dot{\theta}^{m},\dot{\phi}^{m}\rangle d\tau|
\leq (1+\frac{\gamma}{\alpha})\lambda\int_{0}^{t}\|\dot{\phi}^2\|^2 d\tau
     +(1+\frac{\gamma}{\alpha})\lambda\int_{0}^{t}\|\dot{\theta}^{m}\|^2d\tau.
\end{align*}
Combine the ones above, we have
\begin{align*}
 \|\dot{\mathbf{V}}^{m}\|^{2}_{\mathbf{Y}}
\leq -M_{6}\int_{0}^{t}\|\alpha(\tau)\dot{\mathbf{V}}^{m}\|^{2}_{\mathbf{Y}_{\frac{1}{2}}}d\tau+g(\mathbf{V^m}(0)),
\end{align*}
where 
$$
M_{6}=2\min\{\frac{PrC}{Da},1,\frac{1}{\alpha}\}
$$ 
and $\alpha(\tau)=\max\{\frac{Pr}{Da},~(1+\frac{\gamma}{\alpha})\lambda+\frac{Pr Ra^2}{3Da},~(1+\frac{\gamma}{\alpha})\alpha\}+4cM_{so}\|\nabla\theta^m\|$ is nonnegative and local integrable. Thus, by \textbf{Gronwall inequality}, we know
$$
\|\dot{\mathbf{V}}^{m}\|^{2}_{\mathbf{Y}}<\infty, \ \ \ \ \int_{0}^{t}\|\dot{\mathbf{V}}^{m}\|^2_{\mathbf{Y}_{\frac{1}{2}}}d\tau<+\infty,
$$
which means
$\{\mathbf{V}^{m}\}\subset W^{1,2}((0,T),\mathbf{Y}_{\frac{1}{2}})\cap W^{1,\infty}((0,T),\mathbf{Y})$ is bounded. 
According to \textbf{functional analysis}, we obtain that
$$\mathbf{V}^{m}\stackrel{\ast}{\rightharpoonup}\mathbf{V}_{0} \in W^{1,2}((0,T),\mathbf{Y}_{\frac{1}{2}})\cap W^{1,\infty}((0,T),\mathbf{Y}),$$
by the uniqueness of weak limit, we know $\mathbf{V}_{0}$ is the global solution of the equation\eqref{b4}.

%Now, when initial value $\mathbf{V}(0)\in \mathbf{Y}_{1}$, we can choose a sequence $\{\mathbf{V}^{m}(0)\}\in \mathbf{Y}_{\frac{3}{2}}$ such that
%$$\mathbf{V}^{m}(0)\rightarrow \mathbf{V}(0)\in \mathbf{Y}_{1}.$$
%Therefore, we can get a solution sequence
%$$
%\{\mathbf{V}(t,\mathbf{V}^{m}(0))\} \stackrel{\ast}{\rightharpoonup} \mathbf{V}(t,\mathbf{V}(0))\in W^{1,2}((0,T),\mathbf{Y}_{\frac{1}{2}}) \cap W^{1,\infty}((0,T),\mathbf{Y}),
%$$
%because $\{\mathbf{V}(t,\mathbf{V}^{m}(0))\}$ is bounded in $W^{1,2}((0,T),\mathbf{Y}_{\frac{1}{2}})$ $\cap W^{1,\infty}((0,T),\mathbf{Y})$.

The uniqueness of $\mathbf{V}_{0}$ comes from the uniqueness of weak solution according to the theorem \ref{w8}.
\end{proof}

Now, we prove the following corollary which will be used in improving the regularity later.
\begin{corollary}\label{corollary1}The global solution $\mathbf{V}(t)$ of the equation\eqref{b4} belongs to $L^{\infty}((0,T),\mathbf{Y}_{\frac{1}{2}}).$
\end{corollary}
\begin{proof}
After taking the inner product of the equation\eqref{b4} with $\mathbf{V}=(\Delta\psi,\theta,\phi)$ and identical transformations, we have
\begin{align*}
 & \|\nabla\Delta\psi\|^2+\|\nabla\theta\|^2+\|\nabla\phi\|^2
\\
=& -\frac{Da}{PrC}\langle\frac{\partial \Delta\psi}{\partial t},\Delta\psi\rangle
   -\langle\frac{\partial\theta}{\partial t},\theta\rangle-\alpha\langle\frac{\partial\phi}{\partial t},\phi\rangle-\frac{1}{C}\|\Delta\psi\|^2
\\
& -\frac{Ra}{C}\langle\Delta\psi,\frac{\partial \theta}{\partial x}\rangle
  +(\lambda+\gamma\lambda)\langle\phi,\theta\rangle-\lambda\|\theta\|^2-\gamma\lambda\|\phi\|^2.
\end{align*}
After some simple calculations, the following estimated inequalities can be obtained
\begin{align*}
 |\frac{Da}{PrC}\langle\frac{\partial \Delta\psi}{\partial t},\Delta\psi\rangle|
&\leq \frac{1}{2}\|\Delta\psi_{t}\|^2+\frac{Da^2}{2Pr^2C^2}\|\Delta\psi\|^2,
\\
 |\langle\frac{\partial \theta}{\partial t},\theta\rangle|
&\leq \frac{1}{2}\|\theta_{t}\|^2+\frac{1}{2}\|\theta\|^2,
\end{align*}
and 
\begin{align*}
 |\alpha\langle\frac{\partial\phi}{\partial t},\phi\rangle|
&\leq \frac{1}{2}\|\phi_{t}\|^2+\frac{\alpha^2}{2}\|\phi\|^2,
\\
 |\frac{Ra}{C}\langle\Delta\psi,\frac{\partial\theta}{\partial x}\rangle|
&\leq \frac{1}{2}\|\nabla\theta\|^2+\frac{Ra^2}{2C^2}\|\Delta\psi\|^2,
\end{align*}
as well as 
\begin{align*}
 |(\lambda+\gamma\lambda)\langle\phi,\theta\rangle|
&\leq (\frac{\lambda}{4}+\gamma\lambda)\|\phi\|^2+(\lambda+\frac{\gamma\lambda}{4})\|\theta\|^2.
\end{align*}

Combined with the above inequality, whereupon it is concluded that
\begin{align*}
 &\|\nabla\Delta\psi\|^2+\|\nabla\theta\|^2+\|\nabla\phi\|^2
\\
\leq &(\|\Delta\psi_{t}\|^2+\|\theta_{t}\|^2+\|\phi_{t}\|^2) + (\frac{\lambda}{2}+\alpha^2)\|\phi\|^2
\\
 & +( 1+\frac{\gamma\lambda}{2})\|\theta\|^2
   +(\frac{Da^2}{Pr^2C^2}+\frac{2}{C}+\frac{Ra^2}{C^2})\|\nabla\psi\|^2.
\end{align*}
And we have proved that 
$$
\mathbf{V}(t)\in W^{1,\infty}((0,T),\mathbf{Y})\cap W^{1,2}((0,T),\mathbf{Y}_{\frac{1}{2}}),
$$ 
then we can receive 
$$\|\nabla\Delta\psi\|^2+\|\nabla\theta\|^2+\|\nabla\phi\|^2<\infty.$$
This means that
$$\mathbf{V}(t)\in L^{\infty}((0,T),\mathbf{V}_{\frac{1}{2}}).$$
\end{proof}

\section{Existence of attractors in $\mathbf{Y}$}
According to the conclusion in the previous sections, one can see that the equation \eqref{b4} generates a dynamical system $S(t):\mathbf{Y}_{\frac{1}{2}}\rightarrow\mathbf{Y}_{\frac{1}{2}}.$ 
In this section, we aim to prove the existence of attractor of the equation in $\mathbf{Y}$: we firstly prove there exists a bounded absorbing set in $\mathbf{Y}$, and then we demonstrate that semigroup of operators $S(t)$ is uniformly compact in $\mathbf{Y}_{\frac{1}{2}}$. 
%In fact, it is easier to verify the equation \eqref{b4} satisfies the condition (C) in $\mathbf{Y}$. But we need the condition that there exists an absorbing set in $\mathbf{Y}_{\frac{1}{2}}$ to prove the existence of attractor in $\mathbf{Y_{\frac{k}{2}}}(k\geq 2)$.
\begin{theorem}\label{at1}(\textbf{Existence of attractor in $\mathbf{Y}$}) The equation \eqref{b4} has an attractor $\mathcal{A}$ in $\mathbf{Y}$, and $\mathcal{A}$ absorbs all the bounded sets in $\mathbf{Y}$.
\end{theorem}

\begin{proof}
\textbf{Step 1, we will illustrate the existence of bounded absorbing set in $\mathbf{Y}$.}

By considering the second and third equation in the equation \eqref{b4}, we find that we can estimate the two terms $\theta$ and $\phi$ firstly. Thus, after taking the inner product of the second and third equation with $(\theta,\phi)$, we have
\begin{align*}
 \frac{d}{dt}(\frac{1}{2\lambda}\|\theta\|^2+\frac{\alpha}{2\gamma\lambda}\|\phi\|^2)
= -\frac{\|\nabla\theta\|^{2}}{\lambda}-\frac{\|\nabla\phi\|^{2}}{\gamma\lambda}-[\|\theta\|^{2}+\|\phi\|^{2}-2\left\langle \phi,\theta\right\rangle ],
\end{align*}
On the basis of  \textbf{Poincar{\'e} inequality}, we know that
$$
 \|\theta\|^2\leq M_{P}\|\nabla\theta\|^2,
\ \ \ \
 \|\phi\|^2\leq M_{P}\|\nabla\phi\|^2,$$
where $M_{P}>0$ is constant.

Hence, we get the following inequality
\begin{align*}
  \frac{d}{dt}(\frac{1}{2\lambda}\|\theta\|^2+\frac{\alpha}{2\gamma\lambda}\|\phi\|^2)+\widetilde{c} (\frac{1}{2\lambda}\|\nabla\theta\|^2+\frac{\alpha}{2\gamma\lambda}\|\nabla\phi\|^2)\leq -M_{7}(\frac{1}{2\lambda}\|\theta\|^2+\frac{\alpha}{2\gamma\lambda}\|\phi\|^2), 
\end{align*}
where
$$
M_{7}=\frac{1}{M_{P}}\min\{1,\frac{1}{\alpha}\}>0,~~~~0<\widetilde{c}\alpha<1(0<\widetilde{c}<1).
$$

Finally, we take advantage of \textbf{Gronwall inequality} to receive the result 
\begin{align}\label{at2}
\|\theta\|^2+\|\phi\|^2\leq M_{8}e^{-M_{7}t}(\|\theta(0)\|^2+\|\phi(0)\|^2),
\end{align}
and 
\begin{align}\label{at3}
  \int_{0}^{t}\|\nabla\theta\|^2+\|\nabla\phi\|^2 d\tau\leq M_{9}[\|\theta(0)\|^2+\|\phi(0)\|^2],
\end{align}
where 
$$
M_{8}=\frac{\max\{\frac{1}{2\lambda},\frac{\alpha}{2\gamma\lambda}\}}{\min\{\frac{1}{2\lambda},\frac{\alpha}{2\gamma\lambda}\}},~~M_{9}=\frac{\min\{\frac{1}{2\lambda},\frac{\alpha}{2\gamma\lambda}\}}{\min\{\frac{1}{2\lambda},\frac{1}{2\gamma\lambda}\}\widetilde{c}},
$$
which means there exists $t_{0}>0$ such that $\|\theta\|^2+\|\phi\|^2\leq \rho_{0}^2$ ($\rho_{0}$ is a fixed constant) as $t>t_{0}$, i.e. $\|\theta\|$ and $\|\phi\|$ are uniformly bounded.

Similarly, after taking the inner product of the first equation with $\Delta\psi$, we make use of the \textbf{Young inequality} as $t>t_{0}$ to have
\begin{align*}
  \frac{Da}{2Pr}\frac{d}{dt}\|\Delta\psi\|^2
=& -C\|\nabla\Delta\psi\|^2-\|\Delta\psi\|^2
   -Ra\langle\theta,\frac{\partial\Delta\psi}{\partial x}\rangle,
\\
\leq& -\|\Delta\psi\|^2+\frac{Ra^2 }{4C}\|\theta\|^2,
\\
\leq& -\|\Delta\psi\|^2+\frac{Ra^2\rho_{0}^2}{4C}.
\end{align*}
Then, according to the \textbf{Gronwall inequality}, it is concluded that 
\begin{align}\label{at4}
\|\Delta\psi\|^2\leq (1-e^{\frac{2Pr(t_{0}-t)}{Da}})\frac{Ra^2\rho_{0}^2}{4C}.
\end{align}

Thus, there exists $\rho_{R}>0$ such that $\|\Delta\psi\|^2+\|\theta\|^2+\|\phi\|^2<\rho_{R}^{2}$ as $t>t_{0}$, i.e. there exists a bounded absorbing set in $\mathbf{Y}$.
Moreover, by taking the integral
$$
\frac{Da}{2Pr}\frac{d}{dt}\|\Delta\psi\|^2
= -C\|\nabla\Delta\psi\|^2-\|\Delta\psi\|^2
   -Ra\langle\theta,\frac{\partial\Delta\psi}{\partial x}\rangle
$$
from $t_{0}$ to $t$ with respect to t and making use of \eqref{at3}, we have
\begin{align}\label{at5}
\int_{t_{0}}^{t}\|\nabla\Delta\psi\|^2d\tau\leq\frac{Da}{PrC}\|\Delta\psi(t_{0})\|^2+\frac{Ra^2 M_{8}}{C^2M_{7}}[\|\theta(0)\|^2+\|\phi(0)\|^2].
\end{align}

\textbf{Step 2: We will prove $S(t)$ is uniformly compact, i.e. bounded absorbing set in $\mathbf{Y_{\frac{1}{2}}}$.}

After taking the inner product of the equation\eqref{b4} with $(\Delta^2\psi,\Delta\theta,\Delta\phi)$, we have
\begin{align}\label{buchong1}
\begin{aligned}
 &\frac{d}{dt}[\frac{Da}{Pr}\|\nabla\Delta\psi\|^2+\|\nabla\theta\|^2+\alpha\|\nabla\phi\|^2]
\\
= &-2[C\|\Delta^2\psi\|^2+\|\nabla\Delta\psi\|^2+\|\Delta\theta\|^2+\|\Delta\phi\|^2+\lambda\|\nabla\theta\|^2
\\
& +\gamma\lambda\|\nabla\phi\|^2+Ra\langle\frac{\partial\theta}{\partial x},\Delta^2\psi\rangle+\lambda\langle\phi,\Delta\theta\rangle
  +\gamma\lambda\langle\theta,\Delta\phi\rangle+\langle J(\psi,\theta),\Delta \theta\rangle].
\end{aligned}
\end{align}
We use \textbf{Young inequality} to get inequalities
\begin{align}\label{buchong2}
 |Ra\langle\frac{\partial\theta}{\partial x},\Delta^2\psi\rangle|
\leq C\|\Delta^2\psi\|^2 + \frac{Ra^2}{4C}\|\nabla\theta\|^2, 
\end{align}
and 
\begin{align}\label{buchong3}
\lambda|\langle\phi,\Delta\theta\rangle|
\leq \|\nabla\theta\|^2+\frac{\lambda^2}{4}\|\nabla\phi\|^2, 
\end{align}
as well as 
\begin{align}\label{buchong4}
\gamma\lambda|\langle\theta,\Delta\Phi\rangle|
\leq \|\Delta\phi\|^2 + \frac{\gamma^2\lambda^2}{4}\|\theta\|^2.
\end{align}
In the meantime, by \textbf{Sobolev embedding theorem}, we have
\begin{align}\label{buchong5}
  \begin{aligned}
 |\langle J(\psi,\theta),\Delta\theta\rangle|
\leq& \|\nabla\theta\|^2\|\psi\|^2_{H^2} + \|\Delta\theta\|^2,
\\
\leq& cM_{so}\|\nabla\theta\|^2\|\Delta\psi\|^2+\|\Delta\theta\|^2.
  \end{aligned}
\end{align} 

Combined with inequality \eqref{buchong1}-\eqref{buchong5}, the following results are obtained
\begin{align}\label{at6}
  \begin{aligned}
&\frac{d}{dt}(\frac{Da}{Pr}\|\nabla\Delta\psi\|^2+\|\nabla \theta\|^2+\alpha\|\nabla\phi\|^2)
\\
\leq &M_{10}(\frac{Da}{Pr}\|\nabla\Delta\psi\|^2+\|\nabla \theta\|^2+\alpha\|\nabla\phi\|^2)+\frac{(\lambda^2+\gamma^2\lambda^2)\rho_{R}^2}{2},
  \end{aligned}
\end{align}
where 
$$
M_{10}=\frac{Ra^2}{2C}+2\lambda+2+2cM_{so}\|\Delta\psi\|^2+\frac{4\gamma\lambda+\lambda^2}{2\alpha}.
$$
By \textbf{uniform Gronwall inequality}\cite{Temam2013}(page 90), we set
$$
y=\frac{Da}{Pr}\|\nabla\Delta\psi\|^2+\|\nabla \theta\|^2+\alpha\|\nabla\phi\|^2,\ \ \ \ 
g=M_{10},\ \ \ \ 
h=\frac{(\lambda^2+\gamma^2\lambda^2)\rho_{R}^2}{2},
$$
then according to \eqref{at4}-\eqref{at5}, we have
$$
a_{1}=4\int_{t}^{t+r}g(\tau)d\tau\leq M_{11}r.
$$
Formulae \eqref{at3} and \eqref{at5} are further used, we obtain that
$$
a_{3}=\int_{t}^{t+r}y(\tau)d\tau\leq M_{12}r,
$$
where
\begin{align*}
M_{11} &=2cM_{so}\rho_{R}^2+\frac{Ra^2}{2C}+\frac{4\gamma\lambda+\lambda^2}{2\alpha}+2\lambda+2,
\\
M_{12} &=\frac{Da^2}{Pr^2C}\|\Delta\psi(0)\|^2+(\frac{Ra^2M_{8}Da}{C^2M_{7}}Pr+(1+\alpha)M_{9})(\|\theta(0)\|^2+\|\phi(0)\|^2)
\end{align*}
It is easy to know
$$
a_{2} = \int_{t}^{t+r}h(\tau)d\tau =\frac{(\lambda^2+\gamma^2\lambda^2)\rho_{R}^2r}{2} ,
$$
and then we receive that
\begin{align*}
\frac{Da}{Pr}\|\nabla\Delta\psi\|^2+\|\nabla\theta\|^2+\alpha\|\nabla\phi\|^2\leq(\frac{a_{3}}{r}+a_{2})e^{a_{1}},
\end{align*}
which means there exists a bounded absorbing set in $\mathbf{Y_{\frac{1}{2}}}$.

Thus, the proof is completed.
\end{proof}

\section{Existence of $C^{\infty}-$attractor}
In this section, we aim to improve the regularity of the attractor. 
Thus, we need to improve the regularity of the global solution of the equation \eqref{b4} and also prove there exists a bounded absorbing set in $\mathbf{Y_{\frac{k}{2}}}(k\geq 2).$
In order to improve the regularity, we now perform $\Delta^{-1}$ on both sides of the first equation of the equation \eqref{b4} and rewrite it into abstract form as follows:
\begin{align}\label{C1}
\frac{d\mathbf{V}}{dt}=L\mathbf{V}+G(\mathbf{V}),
\end{align}
where
\begin{align*}
  L \mathbf{V}
= \begin{pmatrix}
   \frac{Pr C}{Da} (\Delta^{-1})\Delta^{2} \psi
\\
   \Delta \theta
\\
   \frac{1}{\alpha}\Delta \phi
\end{pmatrix},
\ \ \ \
  G (\mathbf{V})
= \begin{pmatrix}
  -  \frac{Pr }{Da}\psi + \frac{Pr Ra}{Da} (\Delta^{-1}) \frac{\partial \theta}{\partial x}
\\
   \lambda\left(\phi-\theta\right) - J \left(\psi,\theta\right)
\\
  \frac{\gamma \lambda}{\alpha} \left(\theta - \phi \right).
\end{pmatrix},
\ \ \ \ 
\mathbf{V}= (\psi,\theta,\phi).
\end{align*}

\begin{remark}
Because of $\Delta:~H^2(\Omega)\cap H_{0}^{1}\rightarrow L^2(\Omega)$ is isomorphism, $\Delta^{-1}$ is bounded;
\end{remark}
\begin{remark}
Under the condition \eqref{b6}, we can obtain $\Delta^{-1}\Delta\psi=\psi$ and $(\Delta^{-1})\Delta^2\psi$, which means $(\Delta^{-1})\Delta^2\psi$ is an infinitesimal generator of an analytic semigroup.
Thus, operator $L: \mathbf{Y_{1}}\rightarrow \mathbf{Y}$ is infinitesimal generator of the analytic semigroup denoted by $T(t).$
\end{remark}
Based on the fact that the solution $\mathbf{V}\in L^{\infty}((0,T),\mathbf{Y}_{\frac{1}{2}})$, it is easy to know that $G(\mathbf{V})\in L^{1}((0,T),\mathbf{Y})$. Hence, according to the theorem 4.18 in \cite{2011} or theorem in page 259\cite{Pazy1983}, we have the following lemma:
\begin{lemma}\label{C2} The solution of the equation \eqref{C1} $\mathbf{V}(t,\mathbf{V}(0))$ can be expressed as follow:
\begin{align}
\mathbf{V}(t,\mathbf{V}(0))=T(t)\mathbf{V}(0)+\int_{0}^{t}T(t-\tau)G(\mathbf{V}(\tau,\mathbf{V}(0))) d\tau.
\end{align}
\end{lemma}

The lemma \eqref{C2} is very useful for us to improve the regularity of $\mathbf{V}(t).$
And a property of the analytic semigroup $T(t)$ will be used several times as follows:
\begin{align}\label{C3}
\|L^{\beta}T(t)\|\leq M_{\beta}t^{-\beta}e^{-\delta t}~~(\delta>0,~\beta\in \mathbf{R^1}),
\end{align}
where $L^{\beta}$  is the fractional operator generated by $L$ and $M_{\beta}>0$ is constant depending on $\beta$, and the norm $\|\mathbf{V}\|_{\mathbf{Y}_{\beta}}=\|L^{\beta}\mathbf{Y}\|.$

Because the boundary of a rectangular region are almost everywhere smooth, we have the following theorem when assume the boundary of the region studied is smooth.
\begin{theorem}\label{C4}If $\partial\Omega$ is $C^{\infty}$, then the equation \eqref{C1} has a unique solution 
$$
\mathbf{V}(t)\in C^{\infty}((0,T),\mathbf{Y}_{\frac{k}{2}})(k\geq 2)
$$ 
for arbitrary $\mathbf{V}(0)\in \mathbf{X}$. 
\end{theorem}
\begin{proof} We still have two steps to show the results.

\textbf{Step 1: we will improve the regularity of $\mathbf{V}(t)$ with respect to space variables by iteration.}
For the solution $\mathbf{V}(t)\in\mathbf{Y_{\frac{1}{2}}}$, we obviously have
\begin{align*}
 \|J(\psi,\theta)\|^2&\leq 4\int_{\Omega}|\nabla\psi|^2|\nabla\theta|^2dxdz
\\
 &\leq4\|\nabla\psi\|^2_{L^{\infty}(\Omega)}\|\nabla\theta\|^2\leq 4M_{so}[\|\nabla\psi\|^2\|+\|\Delta\psi\|^2]\nabla\theta\|^2,
\end{align*}
which means $G: \mathbf{Y_{\frac{1}{2}}}\rightarrow\mathbf{Y}$ is continuous and bounded. And by \textbf{lemma} \ref{C2},
we can obtain that
\begin{align}
\begin{aligned}
 \|\mathbf{V}(t)\|_{\mathbf{Y_{\beta}}}&\leq\|T(t)\mathbf{V}(0)\|_{\mathbf{Y}_{\beta}}+\int_{0}^{t}\|L^{\beta}T(t-\tau)\|\|G(\mathbf{V})\|d\tau
\\
 &\leq \|T(t)\mathbf{V}(0)\|_{\mathbf{Y}_{\beta}}+M_{\beta}\int_{0}^{t}t^{-\beta}e^{-\delta (t-\tau)} d\tau\|G(\mathbf{V})\|<\infty,
\end{aligned}
\end{align}
for any $\beta\in$ [$\frac{1}{2},1)$.

It reveals that the solution of the equation \eqref{C1}
$$
\mathbf{V}(t)\in L^{\infty}((0,T),\mathbf{Y}_{\beta}).
$$ 
Then there exists a $\theta(\beta)>0$ as $\frac{1}{2}<\beta<1$ such that
\begin{align*}
G:\mathbf{Y}_{\beta}\rightarrow\mathbf{Y}_{\theta(\beta)}~~~~\rm{continuous~~and~~bounded.}
\end{align*}
According to \textbf{lemma} \eqref{C2}, we have
\begin{align*}
\|\mathbf{V}(t)\|_{\mathbf{Y}_{1}}
\leq  \|T(t)\mathbf{V}(0)\|_{\mathbf{Y}_{1}}
     +\int_{0}^{t}\|L^{1-\theta(\beta)}T(t-\tau)\|\|L^{\theta(\beta)}G(\mathbf{V})\|d\tau
<\infty.
\end{align*}
That is to say the solution
$$
\mathbf{V}(t)\in L^{\infty}((0,T),\mathbf{Y}_{\beta})  
\ \
( ~\forall ~\frac{1}{2}\leq\beta\leq 1).
$$

\textbf{Step 2: we will prove the map $ G:\mathbf{Y}_{\frac{k}{2}}\rightarrow\mathbf{Y_{\frac{k-1}{2}}}(k\geq 2)$ is continuous and bounded.}

It is easy to know
\begin{align}\label{C5}
 J(\psi,\theta)
=\frac{\partial \psi}{\partial x}\frac{\partial\theta}{\partial z}
  -\frac{\partial \psi}{\partial z}\frac{\partial\theta}{\partial x}\in H^{k-1}
\end{align}
for $\psi\in H^{k+2}$ and $\theta\in H^{k}$, which means $G:\mathbf{Y}_{\frac{k}{2}}\rightarrow\mathbf{Y_{\frac{k-1}{2}}}(k\geq 2)$ is a continuous and bounded map.

Above all, using an iterative approach, we have
\begin{align}\label{C6}
\mathbf{V}(t)\in L^{\infty}((0,T),\mathbf{Y}_{\frac{k}{2}})(k\geq2).
\end{align}

\textbf{Step 3:  we need to improve the regularity with respect to time variable.} 
%Prove $\mathbf{V}(t)\in C^{\infty}((0,T),\mathbf{Y}_{\frac{k}{2}})(k\geq2)$.
Because of
\begin{align*}
\begin{aligned}
\|\frac{d\mathbf{V}(t)}{dt}\|_{\mathbf{Y}_{\frac{k}{2}}}
&\leq \|LT(t)\mathbf{V}(0)\|_{\mathbf{Y}_{\frac{k}{2}}}+\|G(\mathbf{V})\|_{\mathbf{Y}_{\frac{k}{2}}}+\int_{0}^{t}\|L^{\frac{1}{2}}T(t-\tau)\|\|L^{\frac{k+1}{2}}G(\mathbf{V})\| d\tau
\\
&\leq M_{\frac{k+2}{2}}t^{-\frac{k+2}{2}}e^{-\delta t}\|\mathbf{V}(0)\|+\|G(\mathbf{V})\|_{\mathbf{Y}_{\frac{k}{2}}}+\int_{0}^{t}\|L^{\frac{1}{2}}T(t-\tau)\|\|L^{\frac{k+1}{2}}G(\mathbf{V})\| d\tau
\\
&< \infty 
\end{aligned}
\end{align*}
for $0<t<T$, which means $\mathbf{V}(t)\in C^{1}((0,T),\mathbf{Y}_{\frac{k}{2}})$. 

Thanks to the iterative method again, for any $k\geq 2$, we know the solution
\begin{align*}
\mathbf{V}(t)\in C^{\infty}((0,T),\mathbf{Y}_{\frac{k}{2}}).
\end{align*}

\end{proof}

Finally, we introduce the regularity results of attractors, that is to say the existence of $C^{\infty}-$attractor.
\begin{theorem}
For any $k\geq 1$, there exists an attractor $\mathcal{A}\subset\mathbf{Y}_{\frac{k}{2}}$ absorbs all the bounded sets in $\mathbf{Y}_{\frac{k}{2}}$.
\end{theorem}
\begin{proof}
We have previously  proved that there exists a bounded absorbing set $B\subset{\mathbf{Y_{\frac{1}{2}}}}.$ 
And expression \eqref{C6} combine with lemma \ref{C2} illustrate that there exists a semigroup 
$$
T(t): \mathbf{Y}_{\frac{k}{2}}\rightarrow\mathbf{Y}_{\frac{k}{2}}(k\geq 2).
$$ 
To complete the proof, then we just need to state that there is a bounded absorbing set 
$
\mathcal{B}_{\frac{k}{2}}\subset\mathbf{Y}_{\frac{k}{2}}
$
for $k\geq2$.

For arbitrary initial value $\mathbf{V}(0)\in U_{\frac{k}{2}}\subset\mathbf{Y}_{\frac{k}{2}}$, which $U_{\frac{k}{2}}$ is bounded set in $\mathbf{Y}_{\frac{k}{2}}$, $U_{\frac{k}{2}}\in \mathbf{Y_{\frac{1}{2}}}$ is also bounded. 
Then, there exists $T_{0}>0$ such that for any $T>T_{0}$ we have
\begin{align}\label{C8}
\mathbf{V}\left(T,\mathbf{V}(0)\right)\in B.
\end{align}

Besides, when the time $T$ satisfies $t>T>T_{0}$, we have
\begin{align}\label{C9}
\mathbf{V}(t)=T(t-T)\mathbf{V}(T,\mathbf{V}(0))+\int_{T}^{t}T(t-\tau)G(\mathbf{V})d\tau.
\end{align}
At the same time, we can also get
\begin{align}\label{C10}
\begin{aligned}
  \|T(t-T)\mathbf{V}(T,\mathbf{V}(0))\|_{\mathbf{Y}_{\frac{k}{2}}}
= &\|L^{\frac{k-1}{2}}T(t-T)\mathbf{V}(T,\mathbf{V}(0))\|_{\mathbf{Y}_{\frac{1}{2}}}
\\
\leq & M_{\frac{k-1}{2}}t^{-\frac{k-1}{2}}e^{-\delta t}\|\mathbf{V}\|_{\mathbf{Y}_{\frac{1}{2}}}
\\
= & C_{t}\rightarrow 0 
\end{aligned}
\end{align}
as $t\rightarrow\infty$.
And then we obtain $\|G(\mathbf{V})\|_{\mathbf{Y}_{\frac{k-1}{2}}}<\infty$ by means of the relation \eqref{C5}.

Whereupon, The following inequality is true
\begin{align}\label{C11}
\begin{aligned}
  \|\mathbf{V}\|_{\mathbf{Y}_{\frac{k}{2}}}
\leq & \|T(t-\tau)\mathbf{V}(T,\mathbf{V}(0))\|_{\mathbf{Y}_{\frac{k}{2}}}
      + \int_{T}^{t}\|L^{\frac{1}{2}}T(t-\tau)\|\|G(\mathbf{V})\|_{\mathbf{Y}_{\frac{k-1}{2}}}d\tau
\\
\leq & C_{t}+M_{\frac{1}{2}}\int_{T}^{t}\tau^{-\frac{1}{2}}e^{-\delta \tau} d\tau
       \|G(\mathbf{V})\|_{\mathbf{Y}_{\frac{k-1}{2}}}
\\
<& M_{13}  \ \  \ \  (t\rightarrow\infty),
\end{aligned}
\end{align}
where $M_{13}$ is independent of $\mathbf{V}(0)$.

By the foregoing reason, we have shown there exists a bounded absorbing set $B_{\frac{k}{2}}\subset\mathbf{Y}_{\frac{k}{2}}(k\geq2)$. This proof is completed.

\end{proof}

\section{Conclusions}
In this article, we studied the existence of the $C^\infty$-attractor for a couple stress fluid in saturated porous media. For the target model, we made mathematical deal by the stream function for the 2-dimension incompressible flow.  Afterwards, we obtain the weak solution by Galerkin method and the weak solution is proved to be unique. In order to get the global solution, the regularity of the weak solution is improved with respect to time variable.  So far,  we say the model can generate a dynamic system. Taking advantage of semigroup, the solution can be expressed in the form of integral. In order to prove the existence of C-attractor,  the regularity of the global solution is improved to $C^{\infty}$ by estimating the integral solution and the existence of the bounded absorbing set in higher regularity space is also proved. Finally, by iterative use of uniformly compact condition, the existence of $C^{\infty}$-attractor is proved.  

However, there are two problems worth discussing as follows:

(1) In 2-dimension case, in addition to the advantage of stream function, some key inequalities, which might not be valid in 3-dimension or higher-dimension case, can be proved by Sobolev embedding theorem. Thus, in higher-dimension case, it is much difficult to prove the similar results.

(2) To improve the regularity of solution to $C^{\infty}$, we need the boundary to be smooth and we only assume the boundary smooth.  But, the boundary of rectangular region is almost everywhere smooth.  Thus, it is natural to ask whether for almost everywhere smooth boundary we can improve the regularity to $C^{\infty}$.

In the future work, we will focus on the above two points.

\bibliographystyle{plain}
\bibliography{shidingxing}

\begin{thebibliography}{10}

\bibitem{Straughan2015}
S.~Brian.
\newblock {\em Convection with local thermal non-equilibrium and microfluidic
  effects[M]}.
\newblock Springer-Verlag GmbH, 2015.

\bibitem{Evans2010}
L.~Evans.
\newblock {\em Partial differential equations}.
\newblock American Mathematical Society, Providence, R.I, 2010.

\bibitem{2007An}
S.~N. Gaikwad, M.~S. Malashetty, and K.~R. Prasad.
\newblock An analytical study of linear and non-linear double diffusive
  convection with soret and dufour effects in couple stress fluid.
\newblock {\em International Journal of Non-Linear Mechanics}, 42(7):903--913,
  2007.

\bibitem{Grisvard2011}
P.~Grisvard.
\newblock {\em Elliptic problems in nonsmooth domains}.
\newblock Society for Industrial and Applied Mathematics, jan 2011.

\bibitem{Hill2014}
A.~A. Hill and M.~R. Morad.
\newblock Convective stability of carbon sequestration in anisotropic porous
  media.
\newblock {\em Proceedings of the Royal Society A: Mathematical, Physical and
  Engineering Sciences}, 470(2170):20140373, 2014.

\bibitem{Kumar2021}
P.~Kumar.
\newblock Conditional stability for thermal convection in a rotating
  couple-stress fluid saturating a porous media with temperature- and
  pressure-dependent viscosity using a thermal non-equilibrium model.
\newblock {\em WSEAS Transactions on Heat and Mass Transfer}, 2021.

\bibitem{Luo2012}
H.~Luo.
\newblock Global attractor of atmospheric circulation equations with humidity
  effect.
\newblock {\em Abstract and Applied Analysis}, 2012:1--15, 2012.

\bibitem{Ma2002}
Q.~Ma, S.~Wang, and C.~Zhong.
\newblock Necessary and sufficient conditions for the existence of global
  attractors for semigroups and applications.
\newblock {\em Indiana University Mathematics Journal}, 51(6):1541--1570, 2002.

\bibitem{2007}
T.~Ma.
\newblock {\em Stability and bifurcation of nonlinear evolution equations}.
\newblock science press, Bei jing, 2007.

\bibitem{2011}
T.~Ma.
\newblock {\em Theory and method of partial differential equations}.
\newblock science press, Bei jing, 2011.

\bibitem{MALASHETTY2009781}
M.~S. Malashetty, I.~S. Shivakumara, and S.~Kulkarni.
\newblock The onset of convection in a couple stress fluid saturated porous
  layer using a thermal non-equilibrium model.
\newblock {\em Physics Letters A}, 373(7):781--790, 2009.

\bibitem{Malashetty2009}
M.~S. Malashetty, I.~S. Shivakumara, and S.~Kulkarni.
\newblock The onset of convection in a couple stress fluid saturated porous
  layer using a thermal non-equilibrium model.
\newblock {\em Physics Letters A}, 373(7):781--790, feb 2009.

\bibitem{MINKOWYCZ19993373}
W.~J. Minkowycz, A.~Haji-Sheikh, and K.~Vafai.
\newblock On departure from local thermal equilibrium in porous media due to a
  rapidly changing heat source: the sparrow number.
\newblock {\em International Journal of Heat and Mass Transfer},
  42(18):3373--3385, 1999.

\bibitem{2006Convection}
D.~A. Nield and A.~Bejan.
\newblock {\em Convection in Porous Media, Third Edition}.
\newblock Springer, 2006.

\bibitem{D.A.Nield2017}
D.~A. Nield and A.~Bejan.
\newblock {\em Convection in porous media,5th Edition[M].}
\newblock Springer, 2017.

\bibitem{Nield2002}
D.~A. Nield, A.~V. Kuznetsov, and M.~Xiong.
\newblock Effect of local thermal non-equilibrium on thermally developing
  forced convection in a porous medium.
\newblock {\em International Journal of Heat and Mass Transfer},
  45(25):4949--4955, 2002.

\bibitem{Quangwangjialan2021}
Z.~G. Pan, L.~Jia, Y.~Q. Mao, and Q.~Wang.
\newblock Transitions and bifurcations in couple stress fluid saturated porous
  media using a thermal non-equilibrium model.
\newblock {\em (Sumitted)}.

\bibitem{Pazy1983}
A.~Pazy.
\newblock {\em Semigroups of linear operators and applications to partial
  differential equations}.
\newblock Springer New York, 1983.

\bibitem{Robinson2009}
J.~C. Robinson.
\newblock {\em Infinite-dimensional dynamical systems}.
\newblock Cambridge University Press, May 2009.

\bibitem{Shivakumara2009}
I.~S. Shivakumara.
\newblock Onset of convection in a couple-stress fluid-saturated porous medium:
  effects of non-uniform temperature gradients.
\newblock {\em Archive of Applied Mechanics}, 80(8):949--957, aug 2009.

\bibitem{SONG200953}
L.~Y. Song, Y.~D. Zhang, and T.~Ma.
\newblock Global attractor of the cahn–hilliard equation in hk spaces.
\newblock {\em Journal of Mathematical Analysis and Applications},
  355(1):53--62, 2009.

\bibitem{EliasM.Stein2011}
E.~M. Stein and R.~Shakarchi.
\newblock {\em Functional analysis: introduction to further topics in
  analysis}.
\newblock Princeton University Press, September 2011.

\bibitem{Stokes1984}
V.~K. Stokes.
\newblock Couple stresses in fluids.
\newblock In {\em theories of fluids with microstructure}, pages 34--80.
  Springer Berlin Heidelberg, 1984.

\bibitem{Stokes1984a}
V.~K. Stokes.
\newblock {\em Theories of fluids with microstructure}.
\newblock Springer Berlin Heidelberg, 1984.

\bibitem{Straughan2005}
B.~Straughan.
\newblock Global nonlinear stability in porous convection with a thermal
  non-equilibrium model.
\newblock {\em Proceedings of the Royal Society A: Mathematical, Physical and
  Engineering Sciences}, 462(2066):409--418, dec 2005.

\bibitem{Sunil2014}
Sunil, S.~Choudhary, and A.~Mahajan.
\newblock Conditional stability for thermal convection in a rotating
  couple-stress fluid saturating a porous media with temperature- and
  pressure-dependent viscosity using a thermal non-equilibrium model.
\newblock {\em Journal of Non-Equilibrium Thermodynamics}, 39(2), jan 2014.

\bibitem{Sunil2019}
Sunil, S.~Choudhary, and A.~Mahajan.
\newblock Stability analysis of a couple-stress fluid saturating a porous
  medium with temperature and pressure dependent viscosity using a thermal
  non-equilibrium model.
\newblock {\em Applied Mathematics and Computation}, 340:15--30, 2019.

\bibitem{Sunil2002}
Sunil, R.~C. Sharma, and M.~Pal.
\newblock On a couple-stress fluid heated from below in a porous medium in the
  presence of a magnetic field and rotation.
\newblock {\em Journal of Porous Media}, 5(2):10, 2002.

\bibitem{Temam2013}
R.~Temam.
\newblock {\em Infinite-dimensional dynamical systems in mechanics and
  physics}.
\newblock Springer New York, December 2013.

\bibitem{Vafai2005}
K~Vafai.
\newblock {\em Handbook of porous media}.
\newblock Taylor \& Francis, Boca Raton, 2005.

\bibitem{Zhang_2009}
Y.~Zhang, C.~Zhong, and S.~Wang.
\newblock Attractors in for a class of reaction{\textendash}diffusion
  equations.
\newblock {\em Nonlinear Analysis: Theory, Methods {\&} Applications},
  71(5-6):1901--1908, 2009.

\bibitem{Zhang2008}
Y.~D. Zhang, L.~Y. Song, and T.~Ma.
\newblock The existence of global attractors for 2d navier-stokes equations in
  h k spaces.
\newblock {\em Acta Mathematica Sinica, English Series}, 25(1):51--58, nov
  2008.

\end{thebibliography}

\end{document}